\documentclass[leqno,12pt]{article}

\usepackage{amssymb}
\usepackage{euscript}
\usepackage[dvips]{graphicx}
\usepackage{flafter}
\usepackage{pstricks}
\usepackage{pgf,tikz}
\usetikzlibrary{arrows}

\usepackage{verbatim}

\usepackage{amsmath}

\usepackage{color}

\usepackage{mathrsfs} 

\usepackage{amsthm}

\begin{document}

\baselineskip=18pt
\setcounter{page}{1}

\renewcommand{\theequation}{\thesection.\arabic{equation}}
\newtheorem{theorem}{Theorem}[section]
\newtheorem{lemma}[theorem]{Lemma}
\newtheorem{definition}[theorem]{Definition}
\newtheorem{proposition}[theorem]{Proposition}
\newtheorem{corollary}[theorem]{Corollary}
\newtheorem{fact}[theorem]{Fact}
\newtheorem{problem}[theorem]{Problem}
\newtheorem{conjecture}[theorem]{Conjecture}
\newtheorem{claim}[theorem]{Claim}

\theoremstyle{definition} 
\newtheorem{remark}[theorem]{Remark}
\newtheorem{example}[theorem]{Example}

\newcommand{\eqnsection}{
\renewcommand{\theequation}{\thesection.\arabic{equation}}
    \makeatletter
    \csname  @addtoreset\endcsname{equation}{section}
    \makeatother}
\eqnsection


\def\r{{\mathbb R}}
\def\e{{\mathbb E}}
\def\p{{\mathbb P}}
\def\P{{\bf P}}
\def\E{{\bf E}}
\def\Q{{\bf Q}}
\def\bro{\mathtt{bro}}
\def\z{{\mathbb Z}}
\def\N{{\mathbb N}}
\def\T{{\mathbb T}}
\def\G{G}

\def\ee{\mathrm{e}}



\vglue50pt

\centerline{\large\bf Large deviations for level sets of branching Brownian motion}

\medskip

\centerline{\large\bf and Gaussian free fields}

\bigskip
\bigskip

\centerline{Elie A\"{\i}d\'ekon\footnote{\scriptsize LPMA, Universit\'e Pierre et Marie Curie, 4 place Jussieu, F-75252 Paris Cedex 05, France, {\tt elie.aidekon@upmc.fr} and {\tt zhan.shi@upmc.fr}}, Yueyun Hu\footnote{\scriptsize LAGA, Universit\'e Paris XIII, 99 avenue J-B Cl\'ement, F-93430 Villetaneuse, France, {\tt yueyun@math.univ-paris13.fr}} and Zhan Shi$^{\scriptsize 1}$}

\bigskip
\bigskip

\centerline{Dedicated to the memory of Professor V.N. Sudakov}

\bigskip
\bigskip

{\leftskip=1truecm \rightskip=1truecm \baselineskip=15pt \small

\noindent{\slshape\bfseries Summary.} We study deviation probabilities for the number of high positioned particles in branching Brownian motion, and confirm a conjecture of Derrida and Shi~\cite{bzdev}. We also solve the corresponding problem for the two-dimensional discrete Gaussian free field. Our method relies on an elementary inequality for inhomogeneous Galton--Watson processes.

\bigskip


\bigskip


} 

\bigskip
\bigskip

\section{Introduction}
\label{s:intro}

Consider the model of one-dimensional Branching Brownian Motion (BBM): Initially a particle starts at the origin and performs standard (one-dimensional) Brownian motion. After a random exponential time of parameter 1, the particle splits into two particles; they perform independent Brownian motions. Each of the particles splits into two after an exponential time. We assume that the exponential random variables and the Brownian motions are independent. The system goes on indefinitely. 

Let $X_{\max}(t)$ denote the rightmost position in the BBM at time $t$. McKean~\cite{mckean75} proves that the distribution function of $X_{\max}(t)$ satisfies the F-KPP equation (Fisher~\cite{fisher}, Kolmogorov, Petrovskii and Piskunov~\cite{kolmogorov-petrovski-piskunov}), from which it follows that
$$
\lim_{t\to \infty} \frac{X_{\max}(t)}{t} 
=
2^{1/2},
$$

\noindent in probability. Further order developments can be found in Bramson~\cite{bramson78a} and \cite{bramson83}. For an account of general properties of BBM, see Bovier~\cite{bovier}.
    
The following large deviation estimate for $X_{\max}(t)$ is known (see \cite{rouault}, \cite{chauvin-rouault}): for $x> 2^{1/2}$, 
\begin{equation}
    \lim_{t\to \infty} \, \frac1t \log \p \Big( X_{\max}(t) \ge xt \Big)
    =
    - \psi(x) \, ,
    \label{grandes_dev}
\end{equation}

\noindent where
$$
\psi(x)
:=
\frac{x^2}{2} - 1 \, .
$$

For $x>0$ and $t>0$, let $N(t, x)$ denote the number of particles, in the BBM, alive at time $t$ and positioned in $[tx, \, \infty)$. It is well-known (Biggins~\cite{biggins95}) that for $0<x<2^{1/2}$,
\begin{equation}
    \lim_{t\to \infty} \, \frac{\log N(t,x)}{t}
    =
    1- \frac{x^2}{2},
    \qquad
    \hbox{\rm a.s.}
    \label{biggins}
\end{equation}

\begin{theorem} 
\label{t:bbm}

 Let $x>0$ and $(1- \frac{x^2}{2})^+ < a < 1$. We have 
 $$
 \lim_{t\to \infty} \frac1t \log \p ( N(t, x) \ge \ee^{at}) 
 =
 - I(a,\, x) \, ,
 $$
 
 \noindent where
 $$
 I(a,\, x)
 := 
 \frac{x^2}{2(1-a)} -1\, .
 $$

\end{theorem} 

\medskip


Theorem \ref{t:bbm} gives an affirmative answer to a conjecture by Derrida and Shi~\cite{bzdev}. The conjecture was motivated by a problem for the $N$-BBM, which is a BBM with the additional criterion that the number of particles in the system should never exceed $N$ (whenever the number is more than $N$, the particle at the leftmost position is removed from the system). Let $X^{(N)}_{\max}(t)$ denote the rightmost position in the $N$-BBM at time $t$. It is known (\cite{bzdev}) that 
$$
\psi_N(x)
:=
- \lim_{t\to \infty} \frac1t \log \p( X^{(N)}_{\max}(t) \ge xt),
$$

\noindent exists. In \cite{bzdev}, it is proved that Theorem \ref{t:bbm} implies the following estimate for $\psi_N(x)$:

\medskip

\begin{theorem} 
\label{t:NBBM}

 For $x>2^{1/2}$, we have 
 $$
 \limsup_{N\to \infty} \frac{\log [\psi_N(x) - \psi(x)]}{\log N} 
 \le
 - \Big( \frac{x^2}{2}-1 \Big) \, .
 $$

\end{theorem} 

\medskip

The inequality in Theorem \ref{t:NBBM} is conjectured in \cite{bzdev} to be an equality.

The rest of the paper is as follows. In Section \ref{s:GW}, we present an inequality for inhomogeneous Galton--Watson processes. This inequality will be used in Section \ref{s:pf_bbm} for the proof of Theorem \ref{t:bbm}, and in Section \ref{s:gff} to establish the corresponding result for two-dimensional Gaussian free fields.

\section{An inequality for inhomogeneous Galton--Watson processes}
\label{s:GW}

Let $(Z_n, \, n\ge 0)$ be an inhomogeneous Galton--Watson process, the reproduction law at generation $n$ being denoted by $\nu_n$.\footnote{We write, indifferently, a probability measure $\nu_n$ on $\{ 0, \, 1, \, 2, \ldots \}$ and a random variable whose distribution is $\nu_n$.} More precisely, 
$$
Z_{n+1}
=
\sum_{k=1}^{Z_n} \nu_n^{(k)} ,
\qquad
n\ge 0 \, ,
$$

\noindent where $\nu_n^{(i)}$, $i\ge 1$, are independent copies of $\nu_n$, and are independent of everything up to generation $n$. Let
$$
m_n
:=
\e(\nu_n) .
$$

\noindent We assume $0<m_n<\infty$, for $n\ge 0$.

\medskip

\begin{proposition}
\label{p:GW}

 Let $\alpha>1$ and $n\ge 1$. For all $0\le i<n$, we assume the existence of $\lambda_i >0$ such that 
 \begin{equation}
     \e (\ee^{\lambda_i \nu_i})
     \le
     \ee^{\alpha \lambda_i m_i}.
     \label{hyp} 
 \end{equation}

 \noindent Then for all $\delta>0$ and all integer $\ell \ge 1$,
 \begin{eqnarray*}
  &&\p\Big( Z_n \ge \max\Big\{ \ell, \, (\alpha+\delta)^n \, \ell \max_{0\le i<n} \prod_{j=i}^{n-1} m_j \Big\} \, \Big| \, Z_0 = \ell \Big)
     \\ 
  &\le& n \,  \exp \Big( - \frac{\delta \ell}{\alpha+\delta} \, \min_{0\le i<n} \lambda_i + \max_{0\le i<n} \lambda_i \Big) .
 \end{eqnarray*}

\end{proposition}

\medskip

We say some words about forthcoming applications of the proposition to BBM (in Section \ref{s:pf_bbm}) and to Gaussian free fields (in Section \ref{s:gff}). In both applications, $\alpha+\delta$ is taken to be as close to $1$ as possible, whereas $\ell$ is taken to be $\ee^{\varepsilon n}$ with $\varepsilon>0$ that can be as small as possible (so that $\ell$ is sufficiently large to compensate $\min_{0\le i<n} \lambda_i$ on the right-hand side, but sufficiently small in front of $\max_{0\le i<n} \prod_{j=i}^{n-1} m_j$ on the left-hand side). Roughly speaking, Proposition \ref{p:GW} says that if \eqref{hyp} is satisfied with appropriate $\lambda_i$, then starting at $Z_0 =\ell$, the inhomogeneous Galton--Watson process exceeds $\max\{ \ell, \, \ee^{(1+o(1))n} \max_{0\le i<n} \prod_{j=i}^{n-1} m_j\}$ at generation $n$ with very small probability. 

\bigskip

\noindent {\it Proof of Proposition \ref{p:GW}.} Let $\ell \ge 1$ be an integer. For notational simplification, we write $\p^\ell (\, \cdot \, ) := \p( \, \cdot \, | \, Z_0 = \ell)$. 

Let $b_i \ge 1$, $1\le i\le n$, be integers. We have, for $1\le i< n$,
$$ 
\p^\ell (Z_{i+1} \ge b_{i+1})
\le 
\p^\ell (Z_i \ge b_i ) 
+
\p^\ell\Big( \sum_{k=1}^{b_i} \nu_i^{(k)} \ge b_{i+1} \Big),
$$

\noindent whereas for $i=0$, the inequality simply says $\p^\ell(Z_1 \ge b_1) \le \p^\ell ( \sum_{k=1}^\ell \nu_0^{(k)} \ge b_1)$. By Chebyshev's inequality, 
$$
\p^\ell\Big( \sum_{k=1}^{b_i} \nu_i^{(k)} \ge b_{i+1} \Big)
\le
\ee^{-\lambda_i \, b_{i+1}} [ \e (\ee^{\lambda_i \nu_i})]^{b_i},
$$

\noindent which, by assumption \eqref{hyp}, is bounded by $\exp( - \lambda_i \, b_{i+1} + \alpha b_i \lambda_i m_i)$. Hence
$$ 
\p^\ell (Z_{i+1} \ge b_{i+1})
\le 
\p^\ell (Z_i \ge b_i ) 
+
\exp( - \lambda_i \, b_{i+1} + \alpha \lambda_i m_i b_i).
$$

Let $\delta>0$. We choose $b_0 := \ell$ and, by induction, 
$$
b_{i+1}
:= 
\max\{ \lfloor (\alpha+\delta) m_i b_i \rfloor, \, \ell \}, 
\qquad 
0\le i \le n-1.
$$

\noindent Then $\alpha \lambda_i m_i b_i = \frac{\alpha \lambda_i}{\alpha+\delta} (\alpha+\delta)m_i b_i \le \frac{\alpha \lambda_i}{\alpha+\delta} (1+b_{i+1}) \le \lambda_i + \frac{\alpha \lambda_i}{\alpha+\delta} b_{i+1}$, so
\begin{eqnarray*}
    - \lambda_i \, b_{i+1} + \alpha \lambda_i m_i b_i 
 &\le & - \lambda_i \, b_{i+1} 
    +
    \lambda_i + \frac{\alpha \lambda_i}{\alpha+\delta} b_{i+1}
    \\
 &=& - \frac{\delta \lambda_i}{\alpha+\delta} b_{i+1}
    +
    \lambda_i 
    \\
 &\le&- \frac{\delta \lambda_i}{\alpha+\delta} \ell
    +
    \lambda_i \, . 
\end{eqnarray*}

Consequently, we have, for $1\le i\le n-1$,
$$ 
\p^\ell (Z_{i+1} \ge b_{i+1})
\le 
\p^\ell (Z_i \ge b_i ) 
+
\exp \Big( - \frac{\delta \lambda_i}{\alpha+\delta} \ell + \lambda_i \Big),
$$

\noindent whereas $\p^\ell (Z_1 \ge b_1) \le \ee^{-\lambda_0 \, b_1} [ \e (\ee^{\lambda_0 \nu_0})]^\ell \le \exp( - \frac{\delta \lambda_0}{\alpha+\delta} \ell + \lambda_0)$. Summing over $i$, we obtain:
\begin{eqnarray*}
    \p^\ell (Z_n \ge b_n)
 &\le& \sum_{i=0}^{n-1} \exp \Big( - \frac{\delta \lambda_i}{\alpha+\delta} \ell + \lambda_i \Big)
    \\
 &\le& n \exp \Big( - \frac{\delta \ell}{\alpha+\delta} \min_{0\le i\le n-1} \lambda_i + \max_{0\le i\le n-1} \lambda_i \Big) .
\end{eqnarray*}

\noindent By induction in $n$, $b_n \le \max\{ \ell, \, \max_{0\le i\le n-1} [(\alpha+\delta)^{n-i} \prod_{j=i}^{n-1} m_j] \ell\}$, which is bounded by $\max\{ \ell, \, (\alpha+\delta)^n \ell \max_{0\le i\le n-1} \prod_{j=i}^{n-1} m_j\}$. The proposition follows immediately.\qed

\section{Proof of Theorem \ref{t:bbm}}
\label{s:pf_bbm}

The proof of the theorem relies on the following elementary result, which explains the presence of the constant $I(a, \, x) := \frac{x^2}{2(1-a)} -1$ in the theorem.

\begin{lemma}
\label{l:I(a,x)}

 Let $x>0$ and $(1- \frac{x^2}{2})^+ < a < 1$. We have, for any $t>0$,
\begin{eqnarray}
    \frac1t \sup_{s\in (0, \, t), \; y\le xt : \; (t-s) - \frac{(xt-y)^2}{2(t-s)} = at} \Big( s - \frac{y^2}{2s} \Big) 
 &=& - I(a, \, x) \, .
    \label{I(a,x)_eq1}
    \\
    \frac1t \sup_{s\in (0, \, t), \; y\in \r, \; z\ge xt : \; (t-s) - \frac{(z-y)^2}{2(t-s)} \ge at} \Big( s - \frac{y^2}{2s} \Big) 
 &=& - I(a, \, x) \, .
    \label{I(a,x)_eq2}
\end{eqnarray}

\end{lemma}

\medskip

\noindent {\it Proof.} Clearly, \eqref{I(a,x)_eq2} is a consequence of \eqref{I(a,x)_eq1}: It suffices to observe that for given $(s, \, z)$, the supremum in $y\in \r$ is the supremum in $y\in (-\infty, \, z]$. 

The proof of \eqref{I(a,x)_eq1} is elementary: The maximizer is $s_*= \frac{(1-a)[x^2 -2(1-a)]}{x^2-2(1-a)^2} t$, $y_*= \frac{x}{1-a} s_*$, which is the unique root of the gradient of the Lagrangian, and the supremum is not reached at the boundary.\qed

\bigskip

We often use the elementary Gaussian tail estimate:
$$
\p (|\mathscr{N}| \ge x) 
\le 
\exp \Big( - \frac{x^2}{2\, \mathrm{Var}(\mathscr{N})} \Big),
\qquad
x\ge 0,
$$

\noindent for all mean-zero non-degenerate Gaussian random variable $\mathscr{N}$. As a consequence, for $x\in \r$ and $y\ge 0$,
\begin{equation}
    \p( |\mathscr{N} -x | \le y)
    \le
    \exp \Big( - \frac{x^2}{2\, \mathrm{Var}(\mathscr{N})} + \frac{|x| \, y}{\mathrm{Var}(\mathscr{N})} \Big).
    \label{Gauss}
\end{equation}

\subsection{Lower bound}

The strategy of the lower bound in Theorem \ref{t:bbm} is as follows: Let $\varepsilon>0$. Let $s_*= \frac{(1-a)[x^2 -2(1-a)]}{x^2-2(1-a)^2} t$ and $y_*= \frac{x}{1-a} s_*$ be the maximizer in \eqref{I(a,x)_eq1} of Lemma \ref{l:I(a,x)}. Let the BBM reach $[y_*, \, \infty)$ at time $s_*$ (which, by \eqref{grandes_dev}, happens with probability at least $\exp[- (1+\varepsilon) (\frac{y_*^2}{2s_*} - s_*)] = \ee^{-(1+\varepsilon) I(a, \, x)t}$ for all sufficiently large $t$), then after time $s_*$ the system behaves ``normally" in the sense that by \eqref{biggins}, with probability at least $1-\varepsilon$ for all sufficiently large $t$, the number of descendants positioned in $[xt, \, \infty)$ at time $t$ of the particle positioned in $[y_*, \, \infty)$ at time $s_*$ is at least $\exp\{ (1-\varepsilon) [(t-s_*) - \frac{(xt-y_*)^2}{2(t-s_*)}]\}$ (which is $\ee^{(1-\varepsilon)at}$); note that the condition $0< \frac{xt-y_*}{t-s_*}< 2^{1/2}$ in \eqref{biggins} is automatically satisfied. Consequently, for all sufficiently large $t$,
$$
\p \Big( N(t, x) \ge \ee^{(1-\varepsilon)at}\Big) 
\ge
(1-\varepsilon) \, \ee^{-(1+\varepsilon) I(a, \, x)t} \, . 
$$

\noindent Since $\varepsilon>0$ can be as small as possible, this yields the lower bound in Theorem \ref{t:bbm}.\qed

\subsection{Upper bound}
\label{subs:BBM_ub}

Let $\frac12<\delta<1$. We discretize time by splitting time interval $[0, \, t]$ into intervals of length $t^\delta$: Let $s_i := i t^\delta$ for $0\le i\le M:= t^{1-\delta}$. For notational simplification, we treat $M$ as an integer (upper integer part  should be used for a rigorous treatment; a similar remark applies later when we discretize space).

We first throw away some uninteresting situations. Let $C>0$ be a constant, and let $E_1(t)$ denote the event that all the particles in the BBM lie in $[-Ct, \, Ct]$ at time $s_i$, for all $1\le i\le M$. The expected number of particles that fall out of the interval is bounded by $\sum_{i=1}^M \ee^{s_i} \, \p (\sup_{u\in [0, \, s_i]} |B(u)| \ge Ct)$, where $(B(u), \, u\ge 0)$ denotes a standard one-dimensional Brownian motion. We choose and fix the constant $C>0$ (whose value depends on $a$ and $x$) such that this expected number is $o(\ee^{- I(a, \, x) t})$, $t\to \infty$. By the Markov inequality,
$$
\p(E_1(t)^c)
=
o(\ee^{- I(a, \, x) t}) ,
\qquad
t\to \infty\, .
$$


Let $E_2(t)$ be the event that for all $0\le i< M$, any particle in the BBM alive at time $s_i$ has a total number of descendants fewer than $t^2 \ee^{t^\delta}$ at time $s_{i+1}$. This number has the geometric distribution of parameter $\ee^{-(s_{i+1}-s_i)} = \ee^{-t^\delta}$, i.e., it equals $k$ with probability $(1-\ee^{-t^\delta})^{k-1} \ee^{-t^\delta}$ for all integers $k\ge 1$. By the Markov inequality again, we have
$$
\p(E_2(t)^c)
\le
\sum_{i=0}^{M-1} \ee^{s_i} \sum_{k \ge t^2 \ee^{t^\delta}} (1-\ee^{-t^\delta})^{k-1} \ee^{-t^\delta}
=
o(\ee^{- I(a, \, x) t}) ,
\qquad
t\to \infty\, .
$$

Consequently, for $t\to \infty$,
\begin{equation}
    \p( N(t, x) \ge \ee^{at} )
    \le 
    \p( N(t, x) \ge \ee^{at}, \; E_1(t), \; E_2(t) )
    +
    o(\ee^{- I(a, \, x) t}) \, .
    \label{pf_BBM_ub_eq1}
\end{equation}

We now discretize space. Let $\delta' \in (0, \, \delta)$. [Later, we are going to assume $\delta' < 2\delta-1$.] Let $\varepsilon>0$ be a small constant (which will ultimately go to $0$). Space interval $[-Ct, \, Ct]$ is split into intervals of length $t^{\delta'}$: Let $x_k := k t^{\delta'}$ for $-t^{1-\delta'} \le k\le t^{1-\delta'}$. We call $f: \, \{ s_i, \, 0\le i\le M \} \to \{ x_j := j t^{\delta'}, \, -Ct^{1-\delta'} \le j\le Ct^{1-\delta'}\}$ a path if 
$$
f(0)=0,
\qquad
f(s_M) \ge (x-\varepsilon)t \, .
$$

\noindent The total number of paths is bounded by $(2Ct^{1-\delta'} +1)^{t^{1-\delta}} = \ee^{o(t)}$, $t\to \infty$. 

Consider the BBM. For $1\le i\le M$, a particle at time $s_i$ is said to follow the path $f$ until time $s_i$ if for all $0\le j\le i$, the ancestor of the particle at time $s_j$ lies in $[f(s_j)-t^{\delta'}, \, f(s_j)+t^{\delta'}]$. Let
$$
Z_i(f)
:=
\hbox{\rm number of particles following the path $f$ until time $s_i$} \, .
$$

On the event $E_1(t)$, we have (using the fact that $xt - t^{\delta'} \ge (x-\varepsilon)t$ for all large $t$)
$$
N(t, x) 
\le 
\sum_f Z_M(f)
\le
\# (\mathrm{paths}) \, \max_{f} Z_M(f),
$$

\noindent where $\sum_{f}$ and $\max_{f}$ denote sum and maximum, respectively, over all possible paths $f$, and $\# (\mathrm{paths})$ stands for the total number of paths.

Let $a' \in (0, \, a)$. Since $\# (\mathrm{paths}) = \ee^{o(t)}$ (for $t\to \infty$), it follows that for all sufficiently large $t$ (say $t\ge t_0$), on the event $\{ N(t, x) \ge \ee^{at}\} \cap E_1(t)$, there exists a path $f$ such that $Z_M(f) \ge \ee^{a't}$. Accordingly, for $t\ge t_0$,
\begin{eqnarray*}
    \p( N(t, x) \ge \ee^{at}, \; E_1(t), \; E_2(t))
 &\le& \sum_{f} \p(Z_M(f) \ge \ee^{a't}, \; E_2(t))
    \\
 &\le& \ee^{o(t)} \, \max_{f} \p(Z_M(f) \ge \ee^{a't}, \; E_2(t)) \, .
\end{eqnarray*}

\noindent In view of \eqref{pf_BBM_ub_eq1}, and since $a'$ can be as close to $a$ as possible, the proof of the upper bound in Theorem \ref{t:bbm} is reduced to showing the following: For $x>0$ and $(1- \frac{x^2}{2})^+ < a < 1$,
\begin{equation}
    \limsup_{t\to \infty} \frac1t \max_f \log \p(Z_M(f) \ge \ee^{at}, \; E_2(t))
    \le
    - I(a, \, x) \, ,
    \label{pf_BBM_ub_eq2}
\end{equation}
 
\noindent with $I(a,\, x) := \frac{x^2}{2(1-a)} -1$ as before. [The meaning of $a$ has slightly changed: It is, in fact, $a'$.] 

To bound $\p(Z_M(f) \ge \ee^{at}, \; E_2(t))$, we distinguish two situations. A path $f$ is said to be good if there exists $i\in [1, \, M) \cap \z$ such that
\begin{equation}
    (t-s_i) - \frac{(f(s_M) - f(s_i))^2}{2(t-s_i)} 
    \ge
    (a-\varepsilon)t \, .
    \label{good}
\end{equation}

\noindent It is said to be bad if it is not good. 

When the path $f$ is good, it is easy to bound $\p(Z_M(f) \ge \ee^{at}, \; E_2(t))$; we can even drop $E_2(t)$ in this case: Let $i\in [1, \, M) \cap \z$ be as in \eqref{good}; since $\{ Z_M(f) \ge \ee^{at} \} \subset \{ Z_i(f) \ge 1\}$, we have
\begin{equation}
    \p(Z_M(f) \ge \ee^{at})
    \le
    \e [Z_i(f) ]
    \le
    \ee^{s_i} \, \p(|B(s_i)- f(s_i)| \le t^{\delta'}),
    \label{ub_good_BBM}
\end{equation}

\noindent with $(B(s), \, s\ge 0)$ denoting, as before, a standard Brownian motion. Since $\delta'<\delta<1$, $s_i = i t^\delta$ and $f(s_i) = O(t)$, it follows from \eqref{Gauss} that
$$
\p(|B(s_i)- f(s_i)| \le t^{\delta'})
\le
\exp \Big( - \frac{f(s_i)^2}{2s_i} + o(t) \Big) ,
$$

\noindent uniformly in $i$ and in $f$. This yields that
\begin{eqnarray*}
 &&\p(Z_M(f) \ge \ee^{at}) 
    \\
 &\le& \exp \Big( s_i - \frac{f(s_i)^2}{2s_i} + o(t) \Big) 
    \\
 &\le& \exp \Big\{ \sup_{s\in (0, \, t), \; y\in \r, \; z\ge (x-\varepsilon)t : \; (t-s) - \frac{(z-y)^2}{2(t-s)} \ge (a-\varepsilon)t} \Big( s - \frac{y^2}{2s} \Big) + o(t) \Big\} \, .
\end{eqnarray*}

\noindent By \eqref{I(a,x)_eq2} of Lemma \ref{l:I(a,x)}, the supremum equals $- I(a-\varepsilon, \, x-\varepsilon) t$, as long as $\varepsilon>0$ is sufficiently small such that $x>\varepsilon$ and that $(1- \frac{(x-\varepsilon)^2}{2})^+ < a-\varepsilon$. Hence, uniformly in all good paths $f$,
$$
\limsup_{t\to \infty} \frac1t \log \p(Z_M(f) \ge \ee^{at})
\le
- I(a-\varepsilon, \, x-\varepsilon) \, .
$$

\noindent Since $I(a-\varepsilon, \, x-\varepsilon)$ can be as close to $I(a, \, x)$ as possible, this will settle the case of good paths $f$. To prove \eqref{pf_BBM_ub_eq2}, it suffices to check that, uniformly in all bad paths $f$,
\begin{equation}
    \lim_{t\to \infty} \frac1t \log \p(Z_M(f) \ge \ee^{at}, \; E_2(t))
    =
    - \infty \, .
    \label{pf_BBM_ub_eq3}
\end{equation}

Let $\varepsilon' \in (0, \, \varepsilon)$. For any path $f$, define
$$
\tau = \tau(f, \, t)
:=
\inf\{ i: \, 1\le i\le M, \; Z_i(f) \ge \ee^{\varepsilon' t} \} ,
\qquad
\inf \varnothing := \infty\, .
$$

\noindent On the event $\{ Z_M(f) \ge \ee^{at}\}$, we have $\tau<\infty$, and $Z_\tau(f) \le t^2 \ee^{t^\delta}\ee^{\varepsilon' t}$ on the event $\{ \tau<\infty\} \cap E_2(t)$. Hence
\begin{eqnarray}
 &&\p(Z_M(f) \ge \ee^{at}, \; E_2(t))
    \nonumber
    \\
 &\le&\p(Z_M(f) \ge \ee^{at}, \; Z_\tau(f) \le t^2 \ee^{t^\delta+\varepsilon' t}, \; E_2(t))
    \nonumber
    \\
 &\le& \sum_{i=1}^M \p(Z_M(f) \ge \ee^{at}, \; \ee^{\varepsilon' t} \le Z_i(f) \le t^2 \ee^{t^\delta+\varepsilon' t}, \; E_2(t))
    \nonumber
    \\
 &=& \sum_{i=1}^M \sum_{\ell= \ee^{\varepsilon' t}}^{t^2 \ee^{t^\delta+\varepsilon' t}} \p(Z_M(f) \ge \ee^{at}, \; Z_i(f) = \ell, \; E_2(t)) \, . 
    \label{pf_BBM_ub_eq4}
\end{eqnarray}

Let us have a close look at the probability $\p(Z_M(f) \ge \ee^{at}, \; Z_i(f) = \ell, \; E_2(t))$ on the right-hand side, for $1\le i\le M$ and $\ee^{\varepsilon' t} \le \ell \le t^2 \ee^{t^\delta+\varepsilon' t}$. The sequence $Z_{i+j}(f)$, for $0\le j\le M-i$, can be written as $Z_{j+1} (f) = \sum_{k=1}^{Z_j(f)} \nu_k^{(j)}$, where for each $j$, $\nu_k^{(j)}$, $k\ge 1$, would be i.i.d.\ if the particle at time $s_{i+j}$ were exactly positioned at $f(s_{i+j})$ rather than only lying in the interval $[f(s_{i+j}) - t^{\delta'}, \, f(s_{i+j}) + t^{\delta'}]$. However, $\nu_k^{(j)}$ is stochastically smaller than or equal to $\widetilde{\nu}^{(j)}$, the number of particles in a BBM, starting at position $f(s_{i+j})$, that lie in $[f(s_{i+j+1})- 2t^{\delta'}, \, f(s_{i+j+1}) + 2t^{\delta'}]$ at time $t^\delta$. So we can make a coupling for $(Z_{i+j}(f), \, 0\le j\le M-i)$ and a new process $(\widetilde{Z}_{i+j}(f), \, 0\le j\le M-i)$, which satisfies $\widetilde{Z}_{j+1} (f) = \sum_{k=1}^{\widetilde{Z}_j (f)} \widetilde{\nu}_k^{(j)}$, where for each $j$, $\widetilde{\nu}_k^{(j)}$, $k\ge 1$, are i.i.d.\ having the law of $\widetilde{\nu}^{(j)}$, such that $Z_{i+j}(f) \le \widetilde{Z}_{i+j}(f)$ for all $0\le j\le M-i$. Since $(\widetilde{Z}_{i+j}(f), \, 0\le j\le M-i)$ is an inhomogeneous Galton--Watson process, we can apply Proposition \ref{p:GW}.

Write $\Delta f(s_{i+j}) := f(s_{i+j+1}) - f(s_{i+j})= O(t)$. Note that by \eqref{Gauss},
\begin{eqnarray*}
    \e( \widetilde{\nu}^{(j)} ) 
 &=& \ee^{t^\delta} \, \p ( |B(t^\delta)- \Delta f(s_{i+j})| \le 2t^{\delta'} )
    \\
 &\le& \exp\Big( t^\delta - \frac{(\Delta f(s_{i+j}))^2}{2t^\delta} + O(t^{1+\delta'-\delta}) \Big) 
    \\
 &=:& m_j \, ,  
\end{eqnarray*}

\noindent with $O(t^{1+\delta'-\delta})$ being uniform in $i$, $j$ and $f$. In order to apply Proposition \ref{p:GW}, we need to bound $\max_{0\le k<M-i} \prod_{j=k}^{M-i-1} m_j$, as well as to find a convenient $\lambda_k$ satisfying condition \eqref{hyp} in Proposition \ref{p:GW}. 

Recall that $M:= t^{1-\delta}$. We have, for $0\le k<M-i$, 
\begin{eqnarray*}
    \prod_{j=k}^{M-i-1} m_j 
 &=& \exp\Big( (M-i-k) t^\delta - \frac{1}{2t^\delta} \sum_{j=k}^{M-i-1} (\Delta f(s_{i+j}))^2 + O(t^{2+\delta'-2\delta}) \Big) 
    \\
 &=& \exp\Big( (M-i-k) t^\delta - \frac{1}{2t^\delta} \sum_{j=k}^{M-i-1} (\Delta f(s_{i+j}))^2 + o(t) \Big) ,
\end{eqnarray*}

\noindent as long as $2+\delta'-2\delta <1$ (which is equivalent to $\delta'<2\delta-1$), which we take for granted from now on. By the Cauchy--Schwarz inequality,
$$ 
\sum_{j=k}^{M-i-1} (\Delta f(s_{i+j}))^2 
\ge 
\frac{(f(s_M)- f(s_{i+k}))^2}{M-i-k} \, .
$$

\noindent Recall that $s_j := j t^\delta$ and that $M t^\delta = t$. Hence
\begin{eqnarray*}
    \prod_{j=k}^{M-i-1} m_j 
 &\le& \exp\Big( (M-i-k) t^\delta - \frac{(f(s_M)- f(s_{i+k}))^2}{2(M-i-k)t^\delta} + o(t) \Big) 
    \\
 &=&\exp\Big( (t-s_{i+k}) - \frac{(f(s_M)- f(s_{i+k}))^2}{2(t-s_{i+k})} + o(t) \Big) \, .
\end{eqnarray*}

\noindent If $f$ is a bad path, then by definition of good paths in \eqref{good}, $(t-s_{i+k}) - \frac{(f(s_M)- f(s_{i+k}))^2}{2(t-s_{i+k})} < (a-\varepsilon)t$ for all $k$. Thus 
\begin{equation}
    \max_{0\le k<M-i} \prod_{j=k}^{M-i-1} m_j
    \le
    \ee^{(a-\varepsilon)t + o(t)} \, .
    \label{prod_mj}
\end{equation}

In order to apply Proposition \ref{p:GW}, we still need to find a convenient $\lambda_k$ satisfying condition \eqref{hyp} in the proposition. Let $\alpha>1$. There exists $r>0$ sufficiently small such that $\ee^y \le 1+ \alphaÊy$ for all $y\in [0, \, r]$. On the event $E_2(t)$, we have $\widetilde{\nu}^{(j)} \le t^2 \ee^{t^\delta}$ by definition. Let $\lambda_j := \ee^{-2 t^\delta}$. Then $\lambda_j \widetilde{\nu}^{(j)} \le r$ for all sufficiently large $t$ (and we will be working with such large $t$); hence $\ee^{\lambda_j \widetilde{\nu}^{(j)}} \le 1+ \alpha \lambda_j \widetilde{\nu}^{(j)}$, which yields that
$$
\e(\ee^{\lambda_j \widetilde{\nu}^{(j)}})
\le
1+ \alpha \lambda_j \e(\widetilde{\nu}^{(j)})
\le
1+ \alpha \lambda_j m_j
\le
\ee^{\alpha \lambda_j m_j} \, .
$$

\noindent In words, condition \eqref{hyp} of Proposition \ref{p:GW} is satisfied with the choice of $\lambda_j := \ee^{-2 t^\delta}$. Applying Proposition \ref{p:GW} to $n:= M-i$, we see that for all sufficiently large $t$ and uniformly in $1\le i\le M$ and $\ee^{\varepsilon' t} \le \ell \le t^2 \ee^{t^\delta+\varepsilon' t}$ (recalling that $\varepsilon'<\varepsilon$ and $\varepsilon'<a$) 
$$
\p(Z_M(f) \ge \ee^{at}, \; E_2(t) \, | \, Z_i(f) = \ell)
\le
M \exp ( - c \, \ell \, \ee^{-2 t^\delta} ) \, ,
$$ 

\noindent where $c>0$ is an unimportant constant that does not depend on $t$. A fortiori, $\p(Z_M(f) \ge \ee^{at}, \; E_2(t), \; Z_i(f) = \ell) \le M \exp ( - c \, \ell \, \ee^{-2 t^\delta} )$. By \eqref{pf_BBM_ub_eq4}, we obtain
$$
\p(Z_M(f) \ge \ee^{at}, \; E_2(t))
\le
M^2 \, t^2 \ee^{t^\delta+\varepsilon' t} \exp ( - c \, \ell \, \ee^{-2 t^\delta} ) .
$$

\noindent This yields \eqref{pf_BBM_ub_eq3}, and completes the proof of the upper bound in Theorem \ref{t:bbm}.\qed

\section{Application to discrete Gaussian free fields}
\label{s:gff}

Let $V_N := \{ 1, \ldots, N\}^2$, and $\partial V_N$ be the inner boundary of $V_N$ which is the set of points in $V_N$ having a nearest neighbour outside. Consider the two-dimensional discrete Gaussian free field (GFF) $\Phi = (\Phi(x),\, x\in V_N)$ in $V_N$ with zero boundary conditions as follows: $\Phi$ is a collection of jointly mean-zero Gaussian random variables with $\Phi(x) =0$ for $x\in \partial V_N$ and with covariance given by the discrete Green's function
$$
G_N(x, \, y)
:=
\e_x \Big( \sum_{i=0}^{\tau_{\partial V_N}} {\bf 1}_{\{ S_i =y\} } \Big), 
\qquad
x, \, y\in V_N \backslash \partial V_N\, ,
$$

\noindent where $(S_i, \, i\ge 0)$ is a two-dimensional simple random walk on $\z^2$, $\tau_{\partial V_N}$ the first time the walk hits $\partial V_N$, and $\e_x$ is expectation with respect to $\p_x$ under which $\p_x(S_0=x) =1$. 

In the rest of the paper, we write
\begin{equation}
    \gamma
    :=
    \Big( \frac{2}{\pi} \Big)^{\! 1/2} .
    \label{gamma}
\end{equation}

\noindent This constant originates from the fact that $G_N(0, \, 0) = \gamma^2 \log N + O(1)$, $N\to \infty$ (Lawler \cite{lawler}, Theorem 1.6.6). The maximum of $\Phi$ on $V_N$ was studied by Bolthausen, Deuschel and Giacomin~\cite{bolthausen-deuschel-giacomin}, who proved that
$$
\lim_{N\to \infty} \frac{1}{\log N} \max_{x\in V_N} \Phi(x)
=
2\gamma,
\qquad
\hbox{\rm in probability.}
$$

\noindent [It is possible to have a further development for $\max_{x\in V_N} \Phi(x)$ until constant order of magnitude; see Bramson, Ding and Zeitouni~\cite{bramson-ding-zeitouni}.] Daviaud~\cite{daviaud} was interested in the intermediate level sets
$$
\mathscr{H}_N(\eta)
:=
\{ x\in V_N: \, \Phi(x) \ge 2\gamma \eta \log N \},
\qquad
0<\eta <1 ,
$$

\noindent and proved that for all $0<\eta <1$,
$$
\# \mathscr{H}_N(\eta)
=
N^{2(1-\eta^2) + o(1)},
\qquad
\hbox{\rm in probability,}
$$
 
\noindent where $\# \mathscr{H}_N(\eta)$ denotes the cardinality of $\mathscr{H}_N(\eta)$. Recently, Biskup and Louidor~\cite{biskup-louidor} established the scaling limit of $\mathscr{H}_N(\eta)$ upon an encoding via a point measure.  

We study the deviation probability $\p(\# \mathscr{H}_N(\eta) \ge N^{2a})$, for $1-\eta^2 <a<1$. 

\medskip

\begin{theorem}
\label{t:main}

 Let $\eta\in (0, \, 1)$ and $a\in (1-\eta^2, \, 1)$. We have
 $$
 \p(\# \mathscr{H}_N(\eta) \ge N^{2a})
 =
 N^{- J(a, \, \eta) + o(1)},
 \qquad
 N\to \infty\, ,
 $$
 where
 $$
 J(a, \, \eta)
 :=
 \frac{2\eta^2}{1-a} - 2\, .
 $$
\end{theorem}

\medskip

To prove Theorem \ref{t:main}, let us introduce a useful decomposition. Let $D\subset V_N$ be a square. Define
$$
h_D(x)
:=
\e\big( \Phi(x) \, | \mathscr{F}_{\partial D}) ,
\qquad
x\in D ,
$$

\noindent where $\mathscr{F}_A:= \sigma (\Phi(x), \, x\in A)$ for all $A\subset V_N$, and $\partial D$ denotes the inner boundary of $D$. Let
\begin{equation}
    \Phi^D(x)
    :=
    \Phi(x) - h_D(x), 
    \qquad x \in D .
    \label{PhiD}
\end{equation}
 
\noindent Then $(\Phi^D(x), \, x \in D)$ is independent of $\mathscr{F}_{\partial D \cup D^c}$; in particular, $(\Phi^D(x), \, x \in D)$ and $(h_D(x), \, x \in D)$ are independent. Moreover, $(\Phi^D(x), \, x \in D)$ is a GFF in $D$ in the sense that it is a mean-zero Gaussian field vanishing on $\partial D$ with covariance $\mathrm{Cov} (\Phi^D(x), \, \Phi^D(y)) = \e_x ( \sum_{i=0}^{\tau_{\partial D}} {\bf 1}_{\{ S_i =y\} })$, for $x$, $y\in D\backslash \partial D$, where $\tau_{\partial D}$ is the first hitting time at the inner boundary $\partial D$ by the simple random walk $(S_i)$. 

Write $x_D$ for the centre of $D$. Let
$$
\phi_D 
:= 
\e (\Phi(x_D) \, | \,  \mathscr{F}_{\partial D}) 
=
h_D(x_D).
$$

\noindent [Degenerate case: $\phi_D = \Phi(x)$ if $D = \{ x\}$.] We frequently use an elementary inequality: By Bolthausen, Deuschel and Giacomin~\cite{bolthausen-deuschel-giacomin} p.~1687,
$$
\mathrm{Var} (h_D(x) - \phi_D)
\le
2 \sup_{y\in \partial D} \big| a(x-y)-a(x_D-y)\big|,
$$
 
\noindent where $a(z) := \sum_{n=0}^\infty [ \p_0(S_n=0) - \p_0(S_n=z)]$ with $(S_n, \, n\ge 0)$ denoting as before a simple random walk on $\z^2$. For any $z\in \z^2$, let $|z|$ denote the $L^\infty$-norm of $z$. Since $a(z) = \gamma^2 \log |z| + O(1)$, $|z|\to \infty$ (\cite{lawler}, Theorem 1.6.2), for any $\delta\in (0, \frac14)$,  there exists a constant $c_1 >0$ depending only on $\delta$, such that for all square $D\subset V_N$ with side length $|D|$,
\begin{equation}
    \mathrm{Var} (h_D(x) - \phi_D) \le c_1 \, ,
    \label{BDG_lemma12}
\end{equation}

\noindent uniformly in  $x \in D$ such that $\mathrm{dist}(x, \partial D) \ge \delta \, |D|$ (where $\mathrm{dist}(x, \partial D) := \inf_{y\in \partial D} |x-y|$).

It is possible to estimate $\mathrm{Var} (\phi_D)$. Let $\gamma := (\frac{2}{\pi})^{1/2}$ as in \eqref{gamma}. By equation (7) and Lemma 1 of Bolthausen, Deuschel and Giacomin~\cite{bolthausen-deuschel-giacomin}, there exists a constant $c_2>0$ such that for all square $D\subset V_N$ with side length $m$, 
\begin{equation}
    \mathrm{Var} (\phi_D) 
    \le
    \gamma^2 \log (\frac{N}{m}) + c_2 ,
    \label{Var(phi_D)_ub}
\end{equation}

\noindent and for any $0<\delta<\frac12$, there exists $c_3(\delta)>0$ such that for all square $D\subset V_N$ with $\mathrm{dist}(x_D, \, V_N^c) \ge \delta N$, 
\begin{equation}
    \mathrm{Var} (\phi_D) 
    \ge
    \gamma^2 \log (\frac{N}{m}) - c_3(\delta) .
    \label{Var(phi_D)_lb}
\end{equation}

\noindent [Degenerate case: $m:=1$ if $D$ is a singleton.]

The proof of Theorem \ref{t:main} uses the same ideas as the proof of Theorem \ref{t:bbm} in Section \ref{s:pf_bbm}, with some appropriate modifications. Again, for the sake of clarity, we prove the upper and the lower bounds in distinct paragraphs. The proof is based on the following elementary fact: For $0<\eta<1$ and $1-\eta^2 < a <1$,
\begin{equation}
    \sup_{(s, \, b, \, y): \; 0<s<1, \; y\ge \eta, \; s- \frac{(y-b)^2}{s} \ge a} [(1-s) - \frac{b^2}{1-s}]
    =
    - \Big( \frac{\eta^2}{1-a}-1 \Big) \, .
    \label{prob_variationnel_GFF} 
\end{equation}

\noindent [This is \eqref{I(a,x)_eq2} of Lemma \ref{l:I(a,x)} after a linear transform. The maximizer is $s^* := \frac{a \eta^2}{\eta^2 -(1-a)^2}$, $b^* := \frac{[\eta^2 - (1-a)]\eta}{\eta^2 -(1-a)^2}$, $y^* = \eta$.]

As in the proof for BBM, for notational simplification, we treat several counting quantities (such as $(\log N)^{1-\delta}$ and $N^{1-s_i}$ below) as integers.

\subsection{Upper bound}
\label{subs:ub}

We start with a comparison lemma which implies that  in   $\mathscr{H}_N(\eta) $,  it suffices to  consider those points $x \in V_N$  away from $\partial V_N$. 

\begin{lemma} \label{lemmacomparison}  

Let $\widehat V \subset {\mathbb Z}^2$ be a square containing $V_N$. Let $\widehat \Phi = (\widehat \Phi(x), x \in \widehat V)$ be a Gaussian free field in $\widehat V$ with zero boundary conditions. For any $b \ge 1$, we have $$
\p(\# \mathscr{H}_N(\eta) \ge  b)
\le
2 \p \Big( \# \widehat{\mathscr{H}}_N(\eta) \ge  \frac{b}{2}\Big),$$

\noindent where  
$$ 
\widehat{\mathscr{H}}_N(\eta) :=
\{ x\in V_N: \, \widehat \Phi(x) \ge 2\gamma \eta \log N \}.
$$
\end{lemma}

{\noindent\it Proof of Lemma \ref{lemmacomparison}.} Exactly as in the decomposition in \eqref{PhiD},  there exists a mean-zero Gaussian field $(\widehat h(x), x \in V_N)$ independent of $(\Phi(x), x \in V_N)$,  such that $$ \widehat \Phi(x)= \Phi(x) + \widehat h(x), \qquad x \in V_N.$$

By  symmetry and the independence of $\widehat h(\cdot)$ and $\mathscr{H}_N(\eta)$, \begin{eqnarray*} && \p \Big( \#\big\{ x \in  \mathscr{H}_N(\eta): \widehat h(x) \ge 0 \big\} \ge \frac12 \#  \mathscr{H}_N(\eta)\, \big|  \mathscr{H}_N(\eta) \Big) 
\\
&=&
 \p \Big( \#\big\{ x \in  \mathscr{H}_N(\eta): \widehat h(x) \le 0 \big\} \ge \frac12 \#  \mathscr{H}_N(\eta)\, \big|  \mathscr{H}_N(\eta) \Big).
 \end{eqnarray*}
 
 \noindent Since the sum of the two conditional probabilities is at least  $1$, we  have  $$ \p \Big( \#\big\{ x \in  \mathscr{H}_N(\eta): \widehat h(x) \ge 0 \big\} \ge \frac12 \#  \mathscr{H}_N(\eta)\, \big|  \mathscr{H}_N(\eta) \Big) \ge \frac12.$$

\noindent   Therefore   $\p \big( \# \widehat{\mathscr{H}}_N(\eta) \ge  \frac{b}{2}\big) \ge \p\big( \# \mathscr{H}_N(\eta) \ge  b, \#\big\{ x \in  \mathscr{H}_N(\eta): \widehat h(x) \ge 0 \big\} \ge \frac12 \#  \mathscr{H}_N(\eta)\big) \ge \frac12\,  \p\big( \# \mathscr{H}_N(\eta) \ge  b\big)$, yielding Lemma \ref{lemmacomparison}.\qed

\medskip
In view  of Lemma \ref{lemmacomparison}, to prove the upper bound in Theorem \ref{t:main},   it suffices  to  show that    \begin{equation} \label{newupper}
 \p(\# \mathscr{H}^*_N(\eta) \ge N^{2a})
 \le
 N^{- J(a, \, \eta) + o(1)},
\end{equation}
\noindent where    $ 
\mathscr{H}^*_N(\eta)
=
\{ x\in V^*_N: \,    \Phi(x) \ge 2\gamma \eta \log N \}  $ with \begin{equation}\label{V*N} V^*_N:= \{ x\in V_N: \, \mathrm{dist}(x, \partial V_N) \ge \frac{3N}{8}\}.\end{equation}

Let $\frac56< \delta<1$. Let $L = L(N) := (\log N)^{1-\delta}$. Let $s_0:=1>s_1>...> s_L:=0$ with $s_i-s_{i+1}= (\log N)^{-(1-\delta)}$. 

For $0\le i<L$, let $\mathscr{D}_{s_i}(N)$ denote the partition of $(\frac{N}4)^{2-2s_i}$ squares of side length $(\frac{N}4)^{s_i}$ of $V^*_N$, $\frac{N}{4}$ being the side length of $V^*_N$. [In particular, $\mathscr{D}_{s_0}(N) = \{ V^*_N\}$, the singleton $V^*_N$.] Let $\mathscr{D}_{s_L}(N) := \{ \{ x\}, \, x\in V^*_N\}$, the family of singletons of $V^*_N$. [So for $D= \{ x\} \in \mathscr{D}_{s_L}(N)$, $\phi_D = \Phi(x)$.] We are going to split the family of partitions  $(\mathscr{D}_{s_i}(N))_{i=1}^L$:  Let  $\mathscr{D}^*_{s_0}(N):= \mathscr{D}_{s_0}(N) = \{ V^*_N\}$ and  define recursively for all $1\le i \le L$,  $$\mathscr{D}^*_{s_i}(N)
:= 
\bigcup_{D \in  \mathscr{D}^*_{s_{i-1}}(N)} \{B \in \mathscr{D}_{s_i}(N):   B \subset D, \, \mathrm{dist}(B, \, \partial D) \ge \frac38 \, |D|\},$$

\noindent where $\mathrm{dist}(B,  \, \partial D):= \inf_{x\in B, y \in \partial D} | x- y |$, and $|D|$  denotes as before the side length of $D$.    In particular, $\mathscr{D}^*_{s_L}(N)$ is a collection of singletons.


For any $z \in \z^2$ and $A \subset \z^2$,    let $\Theta_z A:= z +A$  be the set $A$ shifted  by $z$.  Observe  that there exist $z_1, \ldots , z_{4^L} \in \z^2$  (depending on $N, L, (s_i)$ but deterministic)  such that $\max_{1\le j \le 4^L} |z_j| \le \frac{N}{4}$ and that\footnote{As a matter of fact,  \eqref{cover} holds for any $L\ge 1$ and any  $(s_i)_{0\le i \le L}$  and $N$ satisfying  that $s_0=1>s_1>...>s_L=0$ and $N^{s_i- s_{i+1}} \le \frac12$ for all $0\le i < L$.   This fact can be checked by induction on $L$.  To see  the passage from the case $L-1$ to the case  $L$, we denote by $\widetilde V$ the square formed by  aggregating all squares of  $\mathscr{D}^*_{s_1}(N)$.  Clearly  there exist $y_1, .., y_4 \in \z^2$ such that  $\max_{1\le j \le 4} |y_j| \le \frac12 |V^*_N|$ and $V^*_N \subset  \bigcup_{j=1}^4 \Theta_{y_j} \widetilde V$,  we conclude by applying the induction hypothesis to each square  $\Theta_{y_j} \widetilde V$.} \begin{equation}\label{cover}
V^*_N 
\, \subset\, 
\bigcup_{j=1}^{4^L} \Theta_{z_j} \mathscr{D}^*_{s_L}(N) .
\end{equation}

\noindent [By a slight abuse of notation, we have identified, on the right-hand side, $A$ with $\{ x: \, \{x\} \in A\}$ for $A:= \cup_{j=1}^{4^L} \Theta_{z_j} \mathscr{D}^*_{s_L}(N)$.] It follows that $ \p(\# \mathscr{H}_N^*(\eta) \ge N^{2a}) \le \sum_{j=1}^{4^L} \p(\# \mathscr{H}_{L, z_j}(\eta) \ge 4^{-L} N^{2a}   )$, where for any $z \in \z^2$, $$ 
\mathscr{H}_{L, z}(\eta)
:=
\{ x\in \Theta_{z} \mathscr{D}^*_{s_L}(N): \, \Phi(x) \ge 2\gamma \eta \log N \}.$$

  Since $4^L = N^{o(1)}$,  \eqref{newupper} will follow once we prove that
\begin{equation} \label{newupper2}
 \max_{z\in \z^2, |z| \le \frac{N}{4} } \p (\# \mathscr{H}_{L, z}(\eta) \ge N^{2a})
 \le
 N^{- J(a, \, \eta) + o(1)}. 
\end{equation}


Let $C>0$ be a constant. For $z \in \z^2$,   let 
$$
\mathscr{E}_1(N, z)
:=
\Big\{ |\phi_D| \le C \log N, \; \forall 1\le i\le L, \, \forall D \in \Theta_{z}\mathscr{D}^*_{s_i}(N) \Big\} \, , $$

\noindent where $ \Theta_{z}\mathscr{D}^*_{s_i}(N) := \{\Theta_z B: B \in \mathscr{D}^*_{s_i}(N) \}$ for any $1\le i \le L$. This is the analogue for GFF of the event $E_1(t)$ in Section \ref{subs:BBM_ub}.    Since $ \mathrm{Var}   (\phi_D) \le \gamma^2 \log N + c_2$ (see \eqref{Var(phi_D)_ub}) uniformly in $D \in \Theta_{z}\mathscr{D}^*_{s_i}(N)$ and in $z\in \z^2$ such that $|z| \le \frac{N}{4}$,  we can choose $C>0$ sufficiently large such that
$$
\max_{z\in \z^2, \, |z| \le \frac{N}{4}} \p (\mathscr{E}_1(N, z)^c)
=
o(N^{-J(a, \, \eta)}),
\qquad
N\to \infty \, .
$$

Let $\varrho \in (\frac12, \, \frac32-\delta)$. Let
$$
\mathscr{E}_2(N, z)
:=
\Big\{ \max_{B \in \mathtt{ch}(D)} |h_D(x_B) - \phi_D| \le (\log N)^\varrho, \; \forall 1\le i< L, \, \forall D \in \Theta_{z}\mathscr{D}^*_{s_i}(N)\Big\} \, ,
$$

\noindent where $x_B$ denotes as before the centre of the square $B$, and for all $D \in \Theta_{z}\mathscr{D}^*_{s_i}(N)$ with $1\le i<L$,\footnote{For notational simplification, we feel free to omit the dependence of $z$ in $\mathtt{ch}(D)$. The same omission applies to forthcoming quantities such as $Z_i(g)$, ${\bf Z}_i(g)$, $\nu_i^{(D)}$ and $\widetilde{\nu}_i^{(D)}$, without further mention. All probability estimates hold uniformly in $z\in \z^2$ satisfying $|z| \le \frac{N}{4}$.}  
\begin{equation}
    \mathtt{ch}(D)
    =
    \mathtt{ch}(D, \, z)
    :=
    \{ B \subset D \; \mathrm{with} \; B\in \Theta_{z}{\mathscr{D}^*_{s_{i+1}}(N)} \} \, .
    \label{ch(D)}
\end{equation}

\noindent [In words, the elements in $\mathtt{ch}(D)$ play the role of children in the genealogical tree of BBM.] By \eqref{BDG_lemma12} (since $ \mathrm{dist}(x_B, \partial D) \ge  \mathrm{dist}(B, \partial D) \ge \frac{3}{8} |D|$), $\mathrm{Var} (h_D(x_B) - \phi_D) \le c_1$ for all $B\in \mathtt{ch}(D)$ and uniformly in $|z| \le \frac{N}{4}$, which allows to see that
$$
\max_{z\in \z^2, \, |z| \le  \frac{N}{4}} \p (\mathscr{E}_2(N, z)^c)
=
o(N^{-J(a, \, \eta)}),
\qquad
N\to \infty \, .
$$

\noindent Consequently, the following analogue for GFF of \eqref{pf_BBM_ub_eq1} holds: for $N\to \infty$,
\begin{eqnarray}
 && \max_{z\in \z^2, \, |z| \le  \frac{N}{4}} \p (\# \mathscr{H}_{L, z }(\eta) \ge N^{2a})
     \label{pf_GFF_ub_eq1}
    \\
 &\le& \max_{z\in \z^2, \, |z| \le  \frac{N}{4}} \p (\# \mathscr{H}_{L, z  }(\eta) \ge N^{2a}, \; \mathscr{E}_1(N, z), \; \mathscr{E}_2(N, z))
    +
    o(N^{-J(a, \, \eta)})\, . \nonumber
\end{eqnarray}

Let us discretize space. Let $\varepsilon>0$ be a small constant such that $a-\varepsilon>1-(\eta-\varepsilon)^2$. Let $\delta' \in (0, \, \varrho)$. Space interval $[-C\log N, \, C\log N]$ is split into intervals of length $(\log N)^{\delta'}$. We call $g: \, \{ s_i, \, 0\le i\le L \} \to \{ \frac{j}{(\log N)^{1-\delta'}}, \, -C(\log N)^{1-\delta'} \le j\le C(\log N)^{1-\delta'}\}$ a path if 
$$
g(s_0) =0,
\qquad
g(s_L) \ge \eta -\varepsilon \, .
$$

\noindent The total number of paths is $N^{o(1)}$ when $N\to \infty$. 

Define sets of squares ${\bf Z}_0(g) := \Theta_z \{ V^*_N\}$  (the singleton $\{ z+ V^*_N\}$) and for $1\le i \le L$, 
$$
{\bf Z}_i(g)
:=
\Big\{ D \in \Theta_z \mathscr{D}^*_{s_i}(N) : \, |\phi_{D_k} - g(s_k) 2 \gamma \log N| \le 2 \gamma (\log N)^{\delta'}, \; \forall 1\le k \le i \Big\} \, ,
$$

\noindent where $D_k$ denotes the unique square in $\Theta_{z}\mathscr{D}^*_{s_k}(N)$ containing $D$ (so $D_i=D$ for $D \in \Theta_{z}\mathscr{D}^*_{s_i}(N)$). We write
$$
Z_i(g)
:=
\# {\bf Z}_i(g),
\qquad
0\le i\le L \, ,
$$

\noindent the cardinality of ${\bf Z}_i(g)$. On $\mathscr{E}_1(N, z)$, we have $\# \mathscr{H}_{L, z }(\eta) \le \sum_g Z_L(g)$, where $\sum_g$ sums over all possible paths $g$. 

Let $a' \in (0, \, a)$. For all sufficiently large $N$,
\begin{eqnarray*}
 && \max_{z\in \z^2, \, |z| \le \frac{N}{4}} \p ( \# \mathscr{H}_{L, z }(\eta) \ge N^{2a}, \; \mathscr{E}_1(N, z), \; \mathscr{E}_2(N, z))
    \\
 &\le&\# (\mathrm{paths}) \, \max_g \max_{z\in \z^2, |z| \le \frac{N}{4}} \p (Z_L(g) \ge N^{2a'}, \; \mathscr{E}_2(N, z)) \, ,
\end{eqnarray*}

\noindent where $\max_{g}$ denotes maximum over all possible paths $g$, and $\# (\mathrm{paths})$ stands for the total number of paths, which is $N^{o(1)}$ when $N\to \infty$. In view of \eqref{pf_GFF_ub_eq1}, the proof of the upper bound in Theorem \ref{t:main} is reduced to showing the following: For $0<\eta<1$ and $1- \eta^2 < a < 1$,
\begin{eqnarray}
 &&\limsup_{N\to \infty} \frac{1}{\log N} \max_g \max_{z\in \z^2, \, |z| \le \frac{N}{4}} \log \p (Z_L(g) \ge N^{2a}, \; \mathscr{E}_2(N, z))
    \nonumber
    \\
 && \qquad 
    \le - J(a, \, \eta) \, ,
    \label{pf_GFF_ub_eq2}
\end{eqnarray}
 
\noindent with $J(a, \, \eta) := \frac{2\eta^2}{1-a} -2$ as before. 

A path $g$ is said to be good if there exists $i\in [1, \, L) \cap \z$ such that
\begin{equation}
    s_i - \frac{[g(s_i)-g(s_L)]^2}{s_i} 
    \ge
    a-\varepsilon \, .
    \label{good_GFF}
\end{equation}

\noindent [Since $a-\varepsilon>1-(\eta-\varepsilon)^2$, it is clear that $g(s_i)\not= 0$ in this case.] The path is said to be bad if it is not good. 

Let $g$ be a good path. Let $i\in [1, \, L) \cap \z$ be as in \eqref{good_GFF}. We have the following analogue for GFF of \eqref{ub_good_BBM}:
$$
  \p (Z_L(g) \ge N^{2a})
\le
\sum_{D \in  \Theta_{z} \mathscr{D}^*_{s_i}(N)} \, \p \Big\{ |\phi_D| \ge |g(s_i)| 2 \gamma \log N - 2\gamma (\log N)^{\delta'} \Big\} \, .
$$

\noindent Since $g(s_i)\not= 0$, we have $|g(s_i)| 2 \gamma \log N - 2\gamma (\log N)^{\delta'} \ge 0$ by definition of $g$. By \eqref{Var(phi_D)_ub}, $\mathrm{Var}(\phi_D) \le (1-s_i) \gamma^2 \log N+ c_2$ uniformly in   $|z| \le \frac{N}{4}$, so for $D \in \Theta_{z}\mathscr{D}^*_{s_i}(N)$,
\begin{eqnarray*}
 && \p \Big\{ |\phi_D| \ge |g(s_i)| 2 \gamma \log N - 2\gamma (\log N)^{\delta'} \Big\}
    \\
 &\le& \exp\Big( - \frac{(|g(s_i)| 2 \gamma \log N - 2\gamma (\log N)^{\delta'})^2}{2[(1-s_i) \gamma^2 \log N+ c_2]} \Big)
    \\
 &=& \exp\Big( - \frac{2g^2(s_i)}{1-s_i} \log N + o(\log N) \Big) \, ,
\end{eqnarray*}

\noindent uniformly in $i\in [1, \, L) \cap \z$ (recalling that $\frac12<\delta<1$ and that $0<\delta'<\delta$) and in $|z|\le \frac{N}{4}$. Since $\# \mathscr{D}^*_{s_i}(N) \le  N^{2(1-s_i)}$, this yields, uniformly in $|z| \le \frac{N}{4}$,
\begin{eqnarray*}
    \p(Z_L(g) \ge N^{2a})
 &\le& \exp\Big( [2(1-s_i) - \frac{2g^2(s_i)}{1-s_i}] \log N + o(\log N) \Big)
    \\
 &\le& \exp\Big( 2 \sup_{(s, \, b, \, y)} [(1-s) - \frac{b^2}{1-s}] \log N + o(\log N) \Big) ,
\end{eqnarray*}

\noindent the supremum being over $(s, \, b, \, y)$ satisfying $0<s<1$, $y\ge \eta-\varepsilon$ and $s- \frac{(y-b)^2}{s} \ge a-\varepsilon$. By \eqref{prob_variationnel_GFF}, we get that uniformly in good paths $g$, 
$$
\limsup_{N\to \infty} \frac{1}{\log N} \max_{z\in \z^2, \, |z| \le \frac{N}{4}} \log \p (Z_L(g) \ge N^{2a})
\le
- J(a-\varepsilon, \, \eta- \varepsilon) \, .
$$

\noindent As such, the proof of \eqref{pf_GFF_ub_eq2} is reduced to checking that
\begin{equation}
    \lim_{N\to \infty} \frac{1}{\log N} \max_{g \; \mathrm{bad}\; \mathrm{path}} \max_{z\in \z^2, \, |z| \le \frac{N}{4}} \log \p (Z_L(g) \ge N^{2a}, \; \mathscr{E}_2(N, y))
    =
    - \infty \, .
    \label{pf_GFF_ub_eq3}
\end{equation}

Let, for $0\le i<L$ and $D \in \Theta_{z} \mathscr{D}^*_{s_i}(N)$,
\begin{equation}
    \nu_i^{(D)}
:=
\sum_{B \in \mathtt{ch}(D)} {\bf 1}_{\{ |\phi_B - g(s_{i+1}) 2 \gamma \log N| \le 2 \gamma (\log N)^{\delta'}\}} ,
    \label{nu_i_D}
\end{equation}

\noindent where $\mathtt{ch}(D)$ is as in \eqref{ch(D)}. Then 
\begin{equation}
    Z_{i+1}(g)
    = 
    \sum_{D\in {\bf Z}_i(g)} \nu_i^{(D)} ,
    \qquad
    0\le i<L \, .
    \label{branching_structure_prelim}
\end{equation}

\noindent This gives a branching-type process, except that there is lack of independence. So we are going to replace $\nu_i^{(D)}$ by something slightly different. 

Consider two squares $B \subset D$ in $V_N$. Let $\Phi^D(x) := \Phi(x) - h_D(x)$, $x\in D$, as in \eqref{PhiD}. Define
$$
\phi_B^D
:=
 \e  (\Phi^D(x_B)\, | \, \mathscr{F}^D_{\partial B}) \, ,
$$

\noindent where $\mathscr{F}^D_{\partial B} := \sigma (\Phi^D(y), \, y\in \partial B)$, and $x_B$ is as before the centre of $B$. Then $\phi_B^D$ is independent of $h_D(x_B)$, and
\begin{equation}
    \phi_B^D
    =
    \phi_B - h_D(x_B).
    \label{decom-phiBD}
\end{equation}

We now replace $\nu_i^{(D)}$ (defined in \eqref{nu_i_D}) by
$$
\widetilde{\nu}_i^{(D)}
:=
\sum_{B \in \mathtt{ch}(D)} {\bf 1}_{\{ |\phi_B^D  - (g(s_{i+1}) - g(s_i)) 2 \gamma \log N| \le 4\gamma (\log N)^{\delta'} + (\log N)^\varrho \}} .
$$

\noindent Conditionally on ${\bf Z}_i(g)$, the random variables $\widetilde{\nu}_i^{(D)}$, for $D\in {\bf Z}_i(g)$, are independent.

On the event $\mathscr{E}_2(N, y)$, we have $\nu_i^{(D)} \le \widetilde{\nu}_i^{(D)}$, which implies that
$$
Z_i(g)
\le
\widetilde{Z}_i(g),
\qquad
\forall\, 0\le i\le L,
$$

\noindent where $\widetilde{Z}_0(g) := 1$ and for $0\le i<L$,
$$
\widetilde{Z}_{i+1}(g)
= 
\sum_{\ell=1}^{\widetilde{Z}_i(g)} \widetilde{\nu}_i^{(\ell)},
$$

\noindent with $\widetilde{\nu}_i^{(\ell)}$, $\ell \ge 1$, denoting independent copies of $\widetilde{\nu}_i^{(D)}$, which are independent of $\widetilde{Z}_i(g)$. As such, $(\widetilde{Z}_i(g), \, 0\le i\le L)$ is an inhomogeneous Galton--Watson process. 

Let us estimate $\mathrm{Var}(\phi^D_B)$ on the right-hand side. 

Notice that  $ \mathrm{dist}(D, \, \partial V_N) \ge \frac{N}{8}$ and $    \mathrm{dist}(B, \partial D) \ge \frac38 \, |D|$. Recall from \eqref{decom-phiBD} that $\phi_B = \phi_B^D + h_D(x_B)$, the random variables $\phi^D_B$ and $h_D(x_B)$ being independent. So $\mathrm{Var}(\phi^D_B)= \mathrm{Var}(\phi_B) - \mathrm{Var}(h_D(x_B))$. To estimate $\mathrm{Var}(h_D(x_B))$, we use $\mathrm{Var}(X)- \mathrm{Var}(Y) = \mathrm{Var}(X-Y) + 2 \mathrm{Cov}(X-Y, \, Y)$ and $|\mathrm{Cov}(X-Y, Y)|\le [\mathrm{Var}(X-Y)\, \mathrm{Var}(Y)]^{1/2}$, as well as the fact $\mathrm{Var}(h_D(x_B)- \phi_D)\le c_1$ (see \eqref{BDG_lemma12},  since $ \mathrm{dist}(x_B, \partial D) \ge \frac{3}{8} |D|$), to see that
\begin{equation}
    | \mathrm{Var}(h_D(x_B)) - \mathrm{Var}(\phi_D)|
    \le 
    c_1+ 2 \sqrt{c_1\, \mathrm{Var}(\phi_D)}.
    \label{comp-var}
\end{equation}

\noindent Since $\mathrm{Var}(\phi_D) = \gamma^2 (1-s_i) \log N + O(1)$ and $\mathrm{Var}(\phi_B) = \gamma^2 (1-s_{i+1}) \log N + O(1)$ uniformly in $i$, $B$ and $D$ (see \eqref{Var(phi_D)_ub} and \eqref{Var(phi_D)_lb}; this is  where we use the fact that $\mathrm{dist}(D, \, \partial V_N) \ge \frac{N}{8}$), we get
\begin{eqnarray}
    \mathrm{Var}(\phi^D_B)
 &=& \mathrm{Var}(\phi_B) - \mathrm{Var}(h_D(x_B)) 
    \nonumber
    \\
 &=& \gamma^2 (s_i -s_{i+1}) \log N + O((\log N)^{1/2})
    \label{variance(phiDB)1}
    \\
 &=& (1+ O(\frac{1}{(\log N)^{\delta - \frac12}})) \gamma^2 (s_i -s_{i+1}) \log N \, ,
    \label{variance(phiDB)}
\end{eqnarray}

\noindent uniformly in $0\le i<L$ and in   $D \in \Theta_{z}\mathscr{D}_{s_i}(N)$  with  $|z| \le \frac{N}{4}$.

We now estimate $ \e  (\widetilde{\nu}_i)$, where $\widetilde{\nu}_i$ denotes a random variable having the distribution of $\widetilde{\nu}_i^{(D)}$ (for any $D \in \mathscr{D}_{s_i}(N)$). Applying \eqref{Gauss} to $x = (g(s_{i+1}) - g(s_i)) 2 \gamma \log N$ and $y = 4\gamma (\log N)^{\delta'} + (\log N)^\varrho$, and using \eqref{variance(phiDB)} (noting that $\# \mathtt{ch}(D) \le N^{2(s_i-s_{i+1})}$), we arrive at: uniformly in $0\le i<L$ and   in $z\in \z^2$ with  $|z| \le \frac{N}{4}$,
\begin{eqnarray*}
   \e(\widetilde{\nu}_i)
 &\le& N^{2(s_i-s_{i+1})}\, \exp \Big( - \frac{2[g(s_i)-g(s_{i+1})]^2}{s_i-s_{i+1}} \log N 
    \\
 &&\qquad + O((\log N)^{\frac52 - 2\delta}) + O((\log N)^{1+\varrho-\delta}) \Big) .
\end{eqnarray*}

\noindent Note that $(\log N)^{1+\varrho-\delta} = O((\log N)^{\frac52 - 2\delta})$ (because $\varrho< \frac32-\delta$). Recalling $s_L=0$ and $L= (\log N)^{1-\delta}$, we obtain, uniformly in $1\le j < L$ and   in $z\in \z^2$ with  $|z| \le \frac{N}{4}$,
$$
\prod_{i=j}^{L-1}  \e (\widetilde{\nu}_i)
\le
N^{2 s_j - 2 \sum_{i=j}^{L-1} [g(s_i)-g(s_{i+1})]^2/(s_i- s_{i+1}) + o(1)}.
$$

\noindent [This is where the condition $\delta>\frac56$ is needed.] By the Cauchy--Schwarz inequality, $\sum_{i=j}^{L-1} \frac{[g(s_i)-g(s_{i+1})]^2}{s_i- s_{i+1}} \ge \frac{[g(s_j)-g(s_L)]^2}{s_j}$. Since $g$ is a bad path, this yields the following analogue for GFF of \eqref{prod_mj}: uniformly    in     $|z| \le \frac{N}{4}$,
$$
\max_{1\le j<L} \prod_{i=j}^{L-1}  \e (\widetilde{\nu}_i)
\le
N^{2 (a-\varepsilon) + o(1)}.
$$

On the other hand, for any $\varepsilon'\in (0, \, \varepsilon)$ and all sufficiently large $N$,
$$
 \p(\widetilde{Z}_L(g) \ge N^{2a})
\le
\sum_{i=1}^L \sum_{\ell=N^{\varepsilon'}}^{N^{2(s_{i-1}-s_i)+\varepsilon'}}   \p (\widetilde{Z}_L(g) \ge N^{2a}, \; \widetilde{Z}_i(g) = \ell) .
$$

\noindent [This is the analogue for GFF of \eqref{pf_BBM_ub_eq4}.] To apply Proposition \ref{p:GW} to $\p(\widetilde{Z}_L(g) \ge N^{2a} \, | \, \widetilde{Z}_i(g) = \ell)$, we need to find the corresponding $\lambda_j$ (notation of the proposition): Since $\widetilde{\nu}_j \le N^{2(s_j-s_{j+1})} = \ee^{2(\log N)^\delta}$, we can take $\lambda_j := \ee^{-3 (\log N)^\delta}$ (in place of $3$, any constant greater than $2$ will do the job). Applying Proposition \ref{p:GW} to $n:= L-i$, we see that for all sufficiently large $N$ and uniformly in $1\le i\le L$ and $N^{\varepsilon'} \le \ell \le N^{2(s_{i-1}-s_i)+\varepsilon'}$,
$$
\max_{z\in \z^2, \, |z| \le \frac{N}{4}} \p (\widetilde{Z}_L(g) \ge N^{2a}\, | \, \widetilde{Z}_i(g) = \ell)
\le
L \exp ( - c \, \ell \, \ee^{-3 (\log N)^\delta} ) \, ,
$$ 

\noindent where $c>0$ is an unimportant constant. This yields that uniformly    in     $|z| \le \frac{N}{4}$,  $ \p(\widetilde{Z}_L(g) \ge N^{2a} ) \le L^2 N^{2(s_{i-1}-s_i)+\varepsilon'} \exp ( - c \, N^{\varepsilon'} \, \ee^{-3 (\log N)^\delta} )$. Since $Z_L (g) \le \widetilde{Z}_L(g)$ on $\mathscr{E}_2(N, z)$, this yields \eqref{pf_GFF_ub_eq3}, and completes the proof of the upper bound in Theorem \ref{t:main}.\qed 

\subsection{Lower bound}
\label{subs:lb}

Let $0<\eta<1$, $0<b<\eta$, $\varepsilon>0$. Let $0< \zeta< 1$. 

Let $\mathscr{D}_\zeta(N)$ denote the partition of $N^{2-2\zeta}$ squares of side length $N^\zeta$ of $V_N$. For any $D \in \mathscr{D}_\zeta(N)$, let $\widetilde{D} := \{ x \in D:   \mathrm{dist}(x, \partial D) \ge \frac14 N^\zeta \}$ and
\begin{eqnarray*}
    A_D
 &:=& \Big\{ \forall x \in \widetilde D: \, | h_D(x) - \phi_D | \le \varepsilon \log N\Big\} ,
    \\
    B_D
 &:=& \Big\{ \sum_{x\in \widetilde{D}} {\bf 1}_{\{ \Phi^D(x) \ge 2\gamma (\eta- b) \log N \} } \ge N^{2a-\varepsilon} \Big\} .
\end{eqnarray*}

\noindent It is clear that if there exists $D \in \mathscr{D}_\zeta(N)$ such that $\phi_D \ge (2 \gamma b+\varepsilon) \log N$ and that both $A_D$ and $B_D$ are realized, then we have $\# \mathscr{H}_N(\eta) \ge N^{2a-\varepsilon}$. Hence
\begin{eqnarray}
 &&\p\Big( \# \mathscr{H}_N(\eta) \ge N^{2a-\varepsilon} \Big)
    \nonumber
    \\
 &\ge& \p\Big( \exists D \in \mathscr{D}_\zeta(N): \, \phi_D \ge (2 \gamma b+\varepsilon) \log N, \; A_D\cap B_D\Big) . 
    \label{lb:eq1}
\end{eqnarray}

\noindent By Daviaud~\cite{daviaud}, if $\frac{\eta-b}{\zeta} <1$, then for any $D \in \mathscr{D}_\zeta(N)$ and $N\to \infty$,
$$
\p \Big[ \, \#\Big\{ x\in  \widetilde{D}: \, \Phi^D(x) \ge 2 \gamma (\eta-b) \log N \Big\} \ge N^{2\zeta [ 1- \frac{(\eta-b)^2}{\zeta^2} ] -\varepsilon } \Big]
\to
1\, .
$$

\noindent Hence, we have, for all sufficiently large $N$ (say $N\ge N_0$), $\p(B_D) \ge \frac12$ if 
\begin{equation}
    a
    \le
    \zeta \Big( 1- \frac{(\eta-b)^2}{\zeta^2} \Big).
    \label{a}
\end{equation}

\noindent The events $B_D$, $D\in \mathscr{D}_\zeta(N)$, are i.i.d.\ and each $B_D$ is independent of $(\phi_C, \, A_C)$, $C\in \mathscr{D}_\zeta(N)$. We now go back to \eqref{lb:eq1}, and use the fact that 
$$
\p\Big( \bigcup_{i=1}^n (A_i \cap B_i ) \Big) 
\ge
\min_{1\le i \le n} \p(B_i) \,  \p\Big( \bigcup_{j=1}^n A_j \Big),
$$

\noindent if each $B_i$ is independent of $(A_j, \, 1\le j \le n)$. As such, for $N\ge N_0$ (and for $a$ satisfying \eqref{a}),
$$
\p\Big(  \# \mathscr{H}_N(\eta) \ge N^{2a-\varepsilon} \Big)
\ge
\frac12\,  \p\Big( \exists D \in \mathscr{D}_\zeta(N):  \phi_D \ge (2 \gamma b+\varepsilon) \log N,  A_D \Big) \, .
$$

\noindent By \eqref{BDG_lemma12} (since $ \mathrm{dist}(x, \partial D) \ge    \frac{1}{4} |D|$ for any $x \in \widetilde{D}$) and the Gaussian tail, $\p(A_D^c) \le N^{2\zeta} \ee^{ - \varepsilon^2 (\log N)^2/(2 c_1)}$, uniformly in $D \in \mathscr{D}_\zeta(N)$. Hence $\p( \cup_{D \in \mathscr{D}_\zeta(N)} A_D^c) \le N^2 \ee^{ - \varepsilon^2 (\log N)^2/(2 c_1)}$. Consequently, for $a$ satisfying \eqref{a}, any constant $c>0$ and all sufficiently large $N$,
\begin{eqnarray}
 &&\p\Big(  \# \mathscr{H}_N(\eta) \ge N^{2a-\varepsilon} \Big)
    \nonumber
    \\
 &\ge&\frac12\, \p \Big( \exists D \in \mathscr{D}_\zeta(N):\, \phi_D \ge (2 \gamma b+\varepsilon) \log N \Big)
    -
    N^{-c} \, .
    \label{lb:eq2}
\end{eqnarray}

\noindent The probability on the right-hand side is studied in the following lemma.

\medskip

\begin{lemma}
\label{l:lb}

 Let $0 \le \zeta<1$ and $b > 1-\zeta$. Then for $N\to \infty$,
 $$
 \p\Big( \exists D \in \mathscr{D}_\zeta(N):  \phi_D \ge 2\gamma b \log N  \Big)
 =
 N^{2[(1-\zeta)-\frac{b^2}{1-\zeta}] + o(1)} \, . 
 $$

\end{lemma}

\medskip

Admitting Lemma \ref{l:lb} for the moment, we are able to finish the proof of the lower bound in Theorem \ref{t:main}. Indeed, applying Lemma \ref{l:lb} to $b+\frac{\varepsilon}{2\gamma}$ in place of $b$, it follows from \eqref{lb:eq2} that if $b > 1-\zeta$ and $a$ satisfies \eqref{a}, 
$$ 
\p\Big(  \# \mathscr{H}_N(\eta) \ge N^{2a-\varepsilon}\Big)
\ge
N^{2-2\zeta - 2\frac{(b+\frac{\varepsilon}{2\gamma})^2}{1-\zeta} + o(1)},
\qquad
N\to \infty\, .
$$

\noindent The lower bound in Theorem \ref{t:main} follows immediately, with the optimal choice $\eta = \frac{a \eta^2}{\eta^2 -(1-a)^2}$ and $b = \frac{[\eta^2 - (1-a)]\eta}{\eta^2 -(1-a)^2}$.

It remains to prove Lemma \ref{l:lb}. 

\bigskip

\noindent {\it Proof of Lemma \ref{l:lb}.} The argument is quite standard. 

The upper bound, which is not needed in the paper, follows immediately from the Markov inequality, with $\mathrm{Var}(\phi_D)$ being controlled by \eqref{Var(phi_D)_ub}. 

For the lower bound, we only consider those $D$ away from $\partial V_N$: $D \in \mathscr{D}_\zeta(N)$ such that $D \subset V_N^*$ with $V_N^*$ given in \eqref{V*N}. Denoting by $\mathscr{D}^*_\zeta(N)$ the set of such squares $D$. 

Let $K\ge 1$ be a large integer. Define $\zeta_i:= \zeta+ (1-\zeta) \frac{i}{K}$ for $0\le i \le K$. For a square $D \in \mathscr{D}^*_\zeta(N)$, let $D_i$ be the square in $\mathscr{D}^*_{\zeta_i}(N)$ containing $D$ (for $0\le i<K$; so $D_0=D$) and $D_K := V_N^*$.  

Let $$\mathscr{G}_{K, \zeta}(N):= 
\Big\{ D \in \mathscr{D}^*_\zeta(N): \, \forall \, 1\le i \le K, \,   \mathrm{dist}(x_D, \partial D_i) \ge \frac{|D_i|}{4}\Big\},$$

\noindent where, as before, $x_D$ denotes the center of $D$, and $|D_i|$ the side length of $D_i$. Observe that\footnote{To see this, we may use the  construction leading to  \eqref{cover}: let $\mathscr{D}^*_{\zeta_K}(N):= \{V_N^*\}$  and define recursively for $i= N-1, ..., 0$,  $\mathscr{D}^*_{\zeta_i}(N) := \bigcup_{D \in  \mathscr{D}^*_{\zeta_{i+1}}(N)} \{B \in \mathscr{D}_{\zeta_i}(N):   B \subset D, \,   \mathrm{dist}(B, \partial D) \ge \frac38 \, |D|\}$.  Exactly as in \eqref{cover},  $V_N^*$ is covered by the union of at most $4^K$-shifted $\mathscr{D}^*_{\zeta_0}(N)$ squares, so $ \# \mathscr{D}^*_{\zeta_0}(N)\ge 4^{-K} \#  \mathscr{D}^*_\zeta(N)$. The result follows by noting that $\mathscr{D}^*_{\zeta_0}(N) \subset \mathscr{G}_{K, \zeta}(N)$ (let $D \in \mathscr{D}^*_{\zeta_0}(N)$, for any $1\le i \le K $, $ \mathrm{dist}( D_{i-1}, \partial D_i) \ge \frac{3|D_i|}{8}$ by construction, it follows that $ \mathrm{dist}(x_D, \partial D_i)  \ge \frac{3|D_i|}{8} >  \frac{|D_i|}{4}$ as $x_D \in D_{i-1}$). } 
$$ \# \mathscr{G}_{K, \zeta}(N) \ge 4^{-K} \#  \mathscr{D}^*_\zeta(N).$$ 

\noindent Hence $\# \mathscr{G}_{K, \zeta}(N) = N^{2(1-\zeta) +o(1)}$. 

We are going to prove that
\begin{equation}
    \label{proba2} 
    \p\Big( \exists D \in \mathscr{G}_{K, \zeta}(N): \, \phi_D \ge 2\gamma b \log N  \Big)
    \ge
    N^{2[(1-\zeta)-\frac{b^2}{1-\zeta}] + o(1)} \, . 
\end{equation}

For any $1\le i \le K$, we write 
$$ 
\phi_{D_i}
=
c_D(i) \phi_D + Y_D(i),
$$

\noindent where $Y_D(i)$, $1\le i \le K$, is a Gaussian vector independent of $\phi_D$, and
$$
c_D(i)
:=
\frac{\mathrm{Cov}(\phi_{D_i}, \, \phi_D)}{\mathrm{Var}(\phi_D)} \, .
$$

\noindent Since $D\subset D_i$, we can use the decomposition \eqref{decom-phiBD} and in its notation:
$$
\phi_D= \phi^{D_i}_D + h_{D_i}(x_D) \, .
$$

\noindent The independence of $\phi^{D_i}_D$ and $\phi_{D_i}$ gives
\begin{eqnarray}
    \mathrm{Cov}(\phi_{D_i}, \, \phi_D)
 &=& \mathrm{Cov}(\phi_{D_i}, \, h_{D_i}(x_D))
    \nonumber
    \\
 &=& \mathrm{Var}(\phi_{D_i})
    +
    \mathrm{Cov}(\phi_{D_i}, \, h_{D_i}(x_D) - \phi_{D_i}) \, .
    \label{covariance}
\end{eqnarray}

\noindent Let us look at the covariance expression on the right-hand side. By \eqref{Var(phi_D)_ub} and \eqref{Var(phi_D)_lb} (since $D\in \mathscr{G}_{K, \zeta}(N)$), for $0\le i\le K$,
\begin{equation}
    \mathrm{Var}(\phi_{D_i}) = (1-\zeta_i) \gamma^2 \log N + O(1),
    \qquad
    N\to \infty \, ,
    \label{Var(phi(Di))}
\end{equation}

\noindent whereas by \eqref{BDG_lemma12} (since $ \mathrm{dist}(x_D, \, \partial D_i) \ge \frac{|D_i|}{4}$), $\mathrm{Var} (h_{D_i}(x_D) - \phi_{D_i}) \le c_1$. Hence $\mathrm{Cov}(\phi_{D_i}, \, h_{D_i}(x_D) - \phi_{D_i}) = O( (\log N)^{1/2})$ (by Cauchy--Schwarz). Putting this and \eqref{Var(phi(Di))} into \eqref{covariance}, we get
$$
\mathrm{Cov}(\phi_{D_i}, \, \phi_D)
=
(1-\zeta_i) \gamma^2 \log N + O( (\log N)^{1/2}) \, .
$$

\noindent Together with \eqref{Var(phi(Di))} (case $i=0$, so $D_i = D$), this yields
\begin{equation}
    c_D(i)
    =
    \frac{\mathrm{Cov}(\phi_{D_i}, \, \phi_D)}{\mathrm{Var}(\phi_D)}
    =
    \frac{1-\zeta_i}{1-\zeta}
    +
    O( (\log N)^{-1/2}) \, .
    \label{cD(i)}
\end{equation}

Let $\frac12 <\theta < 1$. Let
\begin{eqnarray*}
    I_N
 &:=& [2\gammaÊb \log N , \, 2\gammaÊb \log N + (\log N)^\theta] \, ,
    \\
    \mathscr{A}_D
 &:=& \Big\{ \phi_D \in I_N, \, \max_{1\le i \le K} |Y_D(i)| \le  (\log N)^\theta \Big\}
    \\
 &=& \Big\{ \phi_D \in I_N, \, \max_{1\le i \le K} |\phi_{D_i} - c_D(i) \phi_D| \le  (\log N)^\theta \Big\} \, .
\end{eqnarray*}

\noindent Let
$$
Z
:=
\sum_{D\in \mathscr{G}_{K, \zeta}(N)} {\bf 1}_{\mathscr{A}_D}.
$$

\noindent For each $D\in \mathscr{G}_{K, \zeta}(N)$, $\phi_D$ is independent of $Y_D(i)$, $1\le i\le K$. So
$$
\e(Z)
=
N^{2(1-\zeta)+o(1)} \, \p(\phi_D \in I_N) \, \p\Big( \max_{1\le i \le K} |Y_D(i)| \le (\log N)^\theta\Big).
$$

\noindent By \eqref{Var(phi(Di))} (case $i=0$), $\p(\phi_D \in I_N) = N^{-\frac{2b^2}{1-\zeta} + o(1)}$. On the other hand, $\mathrm{Var}(Y_D(i)) \le \mathrm{Var}(\phi_{D_i}) = O(\log N)$ (by \eqref{Var(phi(Di))}), so $\p ( \max_{1\le i \le K} |Y_D(i)| \le (\log N)^\theta ) \to 1$. It follows that
\begin{equation}
    \e(Z)
    =
    N^{2(1-\zeta) - \frac{2b^2}{1-\zeta} + o(1)} \, .
    \label{E(Z)}
\end{equation}

We now estimate the second moment $\e(Z^2)$.  Observe that 
$$
Z^2
\, 
\le
\, 
Z
+
\sum_{\ell=1}^K 
\sum_{F \in \mathscr{D}^*_{\zeta_\ell}(N)}\; 
\sum_{E,\, E'} \;
\sum_{D, \, D'}
{\bf 1}_{\mathscr{A}_D \cap \mathscr{A}_{D'}},
$$

\noindent where $\sum_{E,\, E'}$ sums over $E$, $E' \in \mathscr{D}^*_{\zeta_{\ell-1}}(N)$ with $E$, $E'\subset F$ and $E\cap E'=\varnothing$, and $\sum_{D, \, D'}$ over $D$, $D'\in \mathscr{D}^*_\zeta(N)$ satisfying $D\subset E$ and $D' \subset E'$ and $ \mathrm{dist}(x_D, \partial E) \ge \frac{|E|}{4}$ and  $ \mathrm{dist}(x_{D'}, \partial E') \ge \frac{|E'|}{4}$. We define
$$
\widetilde{\mathscr{E}}(N)
:=
\Big\{ \forall D \in \mathscr{D}^*_\zeta(N), \, \max_E |h_E(x_D)-\phi_E| \le (\log N)^\theta \Big\} ,
$$

\noindent where $\max_E$ is over all $E \in \mathscr{D}^*_{\zeta_i}(N)$, with $0\le i<K$, such that $E \supset D$ and that $\mathrm{dist}(x_D, \partial E) \ge \frac{|E|}{4}$. The set $\widetilde{\mathscr{E}}(N)$ plays the same role as $\mathscr{E}_2(N,z)$ in the proof of the upper bound. Exactly as for $\mathscr{E}_2(N,z)$, we have, for any constant $c>0$, $\p(\widetilde{\mathscr{E}}(N)^c) = o(N^{- c})$; since $Z^2 \le N^{4(1-\zeta)}$, it follows from \eqref{E(Z)} that
\begin{equation}
    \e(Z^2 \, {\bf 1}_{\widetilde{\mathscr{E}}(N)^c})
    =
    o(\e(Z)),
    \qquad
    N\to \infty\, .
    \label{E(Z2)eq1}
\end{equation}
  
We have
$$
\e(Z^2\, {\bf 1}_{\widetilde{\mathscr{E}}(N)})
\le
\e(Z)
+
\sum_{\ell=1}^K 
\sum_{F \in \mathscr{D}^*_{\zeta_\ell}(N)}\; 
\sum_{E,\, E'} \;
\sum_{D, \, D'}
\p(\mathscr{A}_D \cap \mathscr{A}_{D'} \cap \widetilde{\mathscr{E}}(N)) \, .
$$

\noindent Recall from \eqref{decom-phiBD} that $\phi_D^E = \phi_D - h_E(x_D)$. On the event $\widetilde{\mathscr{E}}(N)$, $h_E(x_D) \le \phi_E + (\log N)^\theta$, so $\phi_D^E \ge \phi_D - \phi_E - (\log N)^\theta$. On the event $\mathscr{A}_D$, $\phi_E \le c_D(\ell-1) \phi_D + (\log N)^\theta$. Consequently, on the event $\mathscr{A}_D \cap \mathscr{A}_{D'} \cap \widetilde{\mathscr{E}}(N)$, we have
\begin{eqnarray*}
    \phi_D^E 
 &\ge& [1- c_D(\ell-1)] \phi_D - 2(\log N)^\theta
    \\
 &\ge& [1- c_D(\ell-1)] 2 \gamma b \log N - 2(\log N)^\theta \, ,
\end{eqnarray*}

\noindent and $\phi_{D'}^{E'} \ge [1- c_{D'}(\ell-1)] 2 \gamma b \log N - 2(\log N)^\theta$ for the same reason. Furthermore, on $\mathscr{A}_D$, 
$$
\phi_F 
\ge 
c_D(\ell) \phi_D - (\log N)^\theta
\ge
c_D(\ell) 2 \gamma b \log N - (\log N)^\theta \, .
$$

\noindent By independence of $\phi_D^E$, $\phi_{D'}^{E'}$ and $\phi_F$, this yields
$$
\p(\mathscr{A}_D \cap \mathscr{A}_{D'} \cap \widetilde{\mathscr{E}}(N)) 
\le
p_{1,N} \, p_{2,N} \, p_{3,N} \, ,
$$

\noindent where
\begin{eqnarray*}
    p_{1,N}
 &:=& \p \{ \phi_D^E \ge [1- c_D(\ell-1)] 2 \gamma b \log N - 2(\log N)^\theta \} ,
    \\
    p_{2,N}
 &:=& \p \{ \phi_{D'}^{E'} \ge [1- c_{D'}(\ell-1)] 2 \gamma b \log N - 2(\log N)^\theta \} ,
    \\
    p_{3,N}
 &:=& \p \{ \phi_F \ge c_D(\ell) 2 \gamma b \log N - (\log N)^\theta \} .
\end{eqnarray*}

\noindent [Note that $p_{1,N}=p_{2,N}$.] By \eqref{variance(phiDB)1} (with $s_i$ and $s_{i+1}$ replaced by $\zeta_{\ell-1}$ and $\zeta$, respectively),
$$
\mathrm{Var}(\phi_D^E)
=
\gamma^2(\zeta_{\ell-1}-\zeta) \log N + O((\log N)^{1/2}) ,
$$

\noindent whereas $\mathrm{Var}(\phi_F) = (1-\zeta_\ell) \gamma^2 \log N + O(1)$ (case $i=\ell$ in \eqref{Var(phi(Di))}), and in view of the value of $c_D(i)$ in \eqref{cD(i)}, we obtain:
\begin{eqnarray*}
    p_{1,N} 
    =
    p_{2,N}
 &\le& \exp\Big( - \frac{2b^2}{1-\zeta} \, \frac{\ell-1}{K}\log N + O((\log N)^\theta) \Big) ,
    \\
    p_{3,N}
 &\le& \exp\Big( - \frac{2b^2}{1-\zeta} \, \frac{K-\ell}{K}\log N + O((\log N)^\theta) \Big) \, .
\end{eqnarray*}

\noindent Consequently,
$$
\e(Z^2\, {\bf 1}_{\widetilde{\mathscr{E}}(N)})
\le
\e(Z)
+
\sum_{\ell=1}^K N^{2(1-\zeta_\ell)} \, N^{4 (\zeta_\ell- \zeta)} \, N^{-\frac{4b^2}{1-\zeta} \, \frac{\ell-1}{K} - \frac{2b^2}{1-\zeta} \, \frac{K-\ell}{K} + o(1)} \, . 
$$

\noindent Note that $2(1-\zeta_\ell) + 4 (\zeta_\ell- \zeta) - \frac{4b^2}{1-\zeta} \, \frac{\ell-1}{K} - \frac{2b^2}{1-\zeta} \, \frac{K-\ell}{K} = 2[(1-\zeta)- \frac{b^2}{1-\zeta}] (1+ \frac{\ell}{K}) + \frac{4b^2}{(1-\zeta)K}$, which is bounded by $2[(1-\zeta)- \frac{b^2}{1-\zeta}] (1+ \frac{1}{K}) + \frac{4b^2}{(1-\zeta)K}$ (for $1\le \ell \le K$; recalling our assumption $b>1-\zeta$ which implies $(1-\zeta)- \frac{b^2}{1-\zeta} <0$). Consequently, for any $\varepsilon>0$, we can choose $K$ sufficiently large such that
$$
\e(Z^2\, {\bf 1}_{\widetilde{\mathscr{E}}(N)})
\le
\e(Z)
+
N^{2[(1-\zeta)- \frac{b^2}{1-\zeta}] + \varepsilon} ,
\qquad
N\to \infty \, .
$$

\noindent Together with \eqref{E(Z2)eq1} and \eqref{E(Z)}, we obtain, for all sufficiently large $N$, $\e(Z^2) \le N^{2\varepsilon} \, \e(Z)$. By the Cauchy--Schwarz inequality,
$$
\p(Z\ge 1)
\ge
\frac{(\e (Z))^2}{\e(Z^2)}
\ge
N^{-2\varepsilon} \, \e(Z) \, .
$$

\noindent In view of \eqref{E(Z)}, this yields the lower bound in Lemma \ref{l:lb}.\qed

\medskip

\begin{remark}

When $\zeta=0$, Lemma \ref{l:lb} gives the following analogue for GFF of \eqref{grandes_dev}: For $b>1$,
$$
 \p\Big( \max_{x\in V_N} \Phi(x) \ge 2\gamma b \log N  \Big)
 =
 N^{2 (1-b^2) + o(1)} ,
 \qquad
 N\to \infty\, . 
 $$

\end{remark}

\bigskip

\noindent {\Large\bf Acknowledgements}

\medskip

We are grateful to Bernard Derrida for stimulating discussions throughout the work. The project was partly supported by ANR MALIN; E.A.\ also acknowledges supports from ANR GRAAL and  ANR Liouville.

\end{document}